\documentclass{article}
\usepackage[arrow,matrix,tips]{xy}
\usepackage{amsmath,amssymb}
\usepackage{amsfonts}
\usepackage{bm}
\def\t-hol{transversely holomorphic}

\def\L{{\cal L}}

\def\bgn{\begin}

\def\CL{\text{\rm CL}}

\def\E{\bold E}
\def\J{{\cal J}}
\def\L{{\cal L}}
\def\Ob{\text{\rm Ob}}
\def\1{{[1]}}
\def\2{{[2]}}
\def\3{{[3]}}
\def\({\left(}
\def\){\right)}
\def\s-circ{\,{\scriptstyle{\circ}}\,}
\def\<<{<\negthinspace \negthinspace<}
\def\Ad{\text{\rm Ad}}

\def\noin{\noindent}

\def\bgn{\begin}

\def\<{<\negthinspace \negthinspace <}
\def\M{\ss M}
\def\metric{generalized metric }

\def\Lam{\Lambda} 
\def\z{\zeta}
\def\sch{\text{\rm Sch}}
\def\rank{\text{\rm rank}}

\title{ \Huge Poisson structures and \\generalized K\"ahler submanifolds\\ 
}

\author{ \Large Ryushi Goto \\ 
Department of Mathematics,
Graduate School of Science,\\ Osaka University, goto@math.sci.osaka-u.ac.jp\\ \\
\it\small Dedicated to Professor Akira Fujiki on the occasion of his sixtieth birthday}

\date{February 22. 2008.}

\begin{document}
\maketitle
\begin{abstract} 
Let $X$ be a compact K\"ahler manifold with a non-trivial holomorphic Poisson structure $\beta$. Then there exist deformations $\{(\J_{\beta t}, \psi_t)\}$  of non-trivial generalized K\"ahler structures with one pure spinor
on $X$.
We prove that every Poisson submanifold of $X$ is a generalized K\"ahler submanifold with respect to 
$(\J_{\beta t}, \psi_t)$ and provide non-trivial examples of generalized K\"ahler submanifolds arising as holomorphic Poisson submanifolds.
We also obtain unobstructed deformations of bi-Hermitian structures constructed from Poisson structures. 
 \end{abstract}

\tableofcontents

\pagenumbering{arabic}

\def\a{\alpha}
\def\b{\beta}
\def\bs{\backslash}
\def\e{\varepsilon}
\def\gam{\gamma}
\def\Gam{\Gamma}
\def\k{\kappa}
\def\del{\delta}
\def\lam{\lambda}
\def\ome{\omega}
\def\Ome{\Omega}
\def\sig{\sigma}
\def\Sig{\Sigma}
\def\eps{\varepsilon}
\def\thet{\theta}
\def\A{{\cal A} }
\def\E{{\cal E}}
\def\f{\bold {f}}
\def\G{G_2}
\def\D{\Bbb D}
\def\Diff{\text{\rm Diff}_0}
\def\Q{\Bbb Q}
\def\R{\Bbb R}
\def\C{\Bbb C}
\def\H{\text{\rm H}}
\def\Z{\Bbb Z}
\def\O{\Bbb O}
\def\P{\frak P}
\def\N{\frak N}
\def\V{{\cal V}}
\def\X{{\cal X}}
\def\w{\wedge}
\def\({\left(}
\def\){\right)}
\def\G{G_2} 
\def\Cl{\frak Cl} 
\def\nab{\nabla}
\def\neg{\negthinspace}
\def\h{\hat}
\def\hrho{\hat{\rho}}
\def\wideh{\widehat}
\def\til{\tilde}
\def\wtil{\widetilde}
\def\ch{\check}
\def\l{\left}
\def\ol{\overline}
\def\pa{\partial}
\def\olpa{\ol{\partial}}
\def\r{\right}
\def\ran{\rangle} 
\def\lan{\langle}
\def\Spin{\text{Spin}(7)}
\def\ss{\scriptscriptstyle}
\def\trian{\triangle}
\def\hook{\hookrightarrow}
\def\hyper{hyperK\"ahler}
\def\Ker{\text{\rm Ker}}
\def\arrow{\longrightarrow}
\def\lrarrow{\Longleftrightarrow}
\def\CP{\dsize{\Bbb C} \bold P}
\def\SL{\ss{\text{\rm SL}}}
\def\Sym{\text{Sym}^3\,\h{\Ome}^1 _X}
\def\tilSym{\widetilde{ \text{Sym}^3}\,\h{\Ome}^1_X }
\def\SymS{\text{Sym}^3\,{\Ome}^1_S }
\def\W{\wedge_{\ome^0}}
\def\sw{{\ss{\w}}}
\def\noin{\noindent}
\def\bsh{\backslash}
\def\reg{\text{\rm reg }}
\def\:{\, :\,}
\def\CL{\text{\rm CL}}
\def\TT{TX\oplus T^*X}
\def\complex{generalized complex }
\def\K\"ahler{generalized K\"ahler}
\newtheorem{theorem}{Theorem}[section]
\newtheorem{definition}[theorem]{Definition}
\newtheorem{lemma}[theorem]{Lemma}
\newtheorem{proposition}[theorem]{Proposition}
\newtheorem{corollary}[theorem]{Corollary}
\newtheorem{example}[theorem]{Example}
\newtheorem{problem}[theorem]{Problem}
\newtheorem{remark}[theorem]{Remark}
\newtheorem{conjecture}[theorem]{Conjecture}
\newenvironment{proof}{{\it Proof}  }{q.e.d.}
\newenvironment{proof*}{ }{ q.e.d}

\large{
\setcounter{section}{-1}
\section{Introduction}
This is a sequel to the paper \cite{Go3} discussing \complex and \K\"ahler manifolds
from a view point of deformations. In section 1.1, we give an exposition of 
\complex structures.
In section 1.2, we introduce a notion of $\J$-submanifolds of a \complex manifold $(X,\J)$ which is due to Oren Ben-Bassart and Mitya Boyarchenko
and  a $\J$-submanifold $M$ inherits the induced generalized complex structure 
$\J_{\M}$ \cite{OM}, \cite{V}.\footnote{Note that there is a different notion of generalized complex submanifolds
due to Gualtieri. To avoid a confusion, we use the terminology of $J$-submanifolds in this paper.}
Both complex submanifolds and symplectic submanifolds arise as special cases of $\J$-submanifolds.}
The notion of $\J$-submanifolds is equivalent under the action of $d$-closed $b$-fields.
We denote by $N^*$ the conormal bundle to $M$ in $X$.
If a submanifold $M$ admits a $\J$-invariant conormal bundle, 
i.e., 
$$\J(N^*)=N^*,$$then $M$ is a $\J$-submanifold.
After a short explanation of \K\"ahler structures in section 1.3,
we prove in section 1.4 that
if a submanifold $M$ of a \K\"ahler manifold $(X,\J_0, \J_1)$ admits a $\J_0$-invariant conormal bundle, then $M$ is also a $\J_1$-submanifold and 
$M$ inherits the induced \K\"ahler structure $(\J_{0,\M}, \J_{1,\M})$
(see theorem 1.9).
In section 1.5, we introduce {\it a \K\"ahler structure with one pure spinor} which is a pair $(\J, \psi)$
consisting of a \complex structure $\J$ and a $d$-closed, non-degenerate, pure spinor $\psi$ such that the induced pair $(\J, \J_\psi)$ is a \K\"ahler structure. 
On a K\"ahler manifold with a K\"ahler form $\ome$, there is the ordinary \K\"ahler structure with one pure spinor, where the pure spinor is given by $\exp(\sqrt{-1}\ome)$.
\\\\
{\bf Theorem 1.14} 
{\it  
Let $(X, \J_0, \psi)$ be a \K\"ahler manifold with one pure spinor.
If a submanifold $M$ admits 
a $\J_0$-invariant conormal bundle, 
then the pull back $i_{\M}^*\psi$ is 
a $d$-closed, non-degenerate, pure spinor on $M$ and the induced pair $(\J_{0,\M}, i_{\M}^*\psi)$ is 
a \K\"ahler structure with one pure spinor on $M$.} 
\\ \\
This is analogous to the fact that the pull back of a K\"ahler form to a complex submanifold is also a K\"ahler form.
In section 1.6, 
let $X$ be a compact K\"ahler manifold with a Poisson structure $\b$ and 
a K\"ahler form $\ome$.
(We always consider holomorphic Poisson structures in this paper.)
A complex submanifold $M$ is a Poisson submanifold if there exists the induced Poisson structure $\b_{\M}$ on $M$ (see Definition 1.17). 
Note that in algebraic geometry, Poisson schemes and Poisson subschemes are developed
\cite {Na}, \cite{Pol}.
A Poisson structure $\b$ on a compact K\"ahler manifold generates 
deformations of \complex structures $\{\J_{\b t}\}$ parametrized by $t$. 
By applying the stability theorem \cite{Go3}, we obtain deformations of \K\"ahler structures $\{ \J_{\b t}, \psi_t\}_{t\in \trian'}$ with one pure spinor, where $\trian'$ is a one dimensional complex disk.
Then it turns out that 
every Poisson submanifold $M$ of $X$ admits a $\J_{\b t}$-invariant conormal bundle, which is 
a \K\"ahler submanifold of $(X, \J_{\b, t}, \psi_t)$ for 
$t\in \trian'$.
Thus every Poisson submanifold inherits the induced \K\"ahler structures 
\\
\\
{\bf Theorem 1.20} {\it  Let $M$ be a Poisson submanifold of a Poisson manifold $X$ with a K\"ahler structure $\ome$. Then $M$ is a \K\"ahler manifold with the induced structure 
$(\J_{\b t,\M}, \psi_{t,\M})$.}
\\
\\
In section 1.7, we discuss examples of \K\"ahler manifolds arising as Poisson submanifolds. 
For instance, 
every hypersurface of degree $d\leq 3$ in $\C P^3$
is a Poisson submanifold which admits a non-trivial \K\"ahler
structure. We also exhibit invariant Poisson submanifolds of Poisson manifolds given by the action of commutative complex group.  
\\
\indent
In section $2$, we will give a formal proof of the stability theorem in the case of 
a K\"ahler manifold with a Poisson structure $\b$.
\footnote{This construction provides a formal family of deformations of \K\"ahler structures 
 in the case of Poisson deformations (see section 4 in \cite{Go3} for the convergence of the formal power series).} 
The construction in the special case is based on 
the Hodge decomposition and the Lefschetz decomposition, 
which is an application of theorem showing unobstructed deformations of 
generalized Calabi-Yau (metrical) structures \cite{Go2}.
(Note that the proof in the general case is similar and depends on the generalized Hodge decomposition.) 
The method of our proof 
is a generalization of the one in unobstructed deformations of Calabi-Yau manifolds due to Bogomolov-Tain-Todorov \cite{Ti}. 
\\
\indent
In section $3$, we discuss an application of the stability theorem to 
deformations of Bihermitian structures. 
By the one to one correspondence between \K\"ahler structures and 
Bihermitian structures with a torsion condition, 
our deformations of \K\"ahler structures $\{\J_{\b t}, \psi_t\}$ gives rise to 
deformations of Bihermitian structures
 $\{I^+(t), I^-(t), h_t\}$. 
In particular, we show that the infinitesimal deformations of 
$\{I^\pm(t)\}$ are respectively given by the class $[\pm\b\cdot \ome]\in H^{0,1}(T^{1,0})$
defined by the contraction of the Poisson structure $\b$ and 
the K\"ahler from $\ome$. In other words,  it follows
from the stability theorem that the class $[\b\cdot\ome]$ gives rise to 
unobstructed deformations of complex structures.
As an example, we discuss deformations of complex structures on the product
of $\C P^1$ and a complex torus and show that Poisson structures generate
the Kuranishi family in terms of the classes $[\b\cdot \ome]$. 
\\
\indent
Hitchin gave a construction of
Bihermitian structures on Del Pezzo surfaces \cite{Hi2}, \cite {Hi3} by 
the Hamiltonian diffeomorphisms. 
Recently Gualtieri extended Hithicn's construction for Poisson manifolds \cite{Gu3}.
It seems to be interesting to compare the construction by the stability theorem and the one by Hamiltonian diffeomorphsims.
It must be noted that the construction by Hamiltonian diffeomorphisms requires that  the Kodaira-Spencer class $[\b\cdot\ome]$ vanishes. 
\footnote{ For instance, on a complex torus and a hyperK\"ahler manifold, there are Poisson structures with 
non-vanishing class $[\b\cdot \ome]$.}
\\
\indent
The author would like to thank Professor Fujiki for valuable comment about 
Bihermitian structures. He would like to express his gratitude 
to Professor Namikawa for his meaningful suggestions about Poisson geometry. 
He also thanks to Professor Vaisman for  a kind comment about definition of 
\complex submanifolds. 
He is grateful to Professor Gualtieri for his comment about deformations of \K\"ahler structures.

\large{ \numberwithin{equation}{section}
\section{Generalized K\"ahler submanifolds}
\subsection{Generalized complex structures}
Let $X$ be a real compact manifold of $2n$ dimension.
We denote by $\pi : \TT \to TX$ the projection to the first component.
The natural coupling between $TX$ and $T^*X$ defines 
the symmetric bilinear form $\lan\, ,\,\ran$ on the direct sum $\TT$.
Then we have the fibre bundle SO$(\TT)$ over $X$ with fibre the special orthogonal group 
with respect to $\TT$.
Note that SO$(\TT)$ is a subbundle of End$(\TT)$.
{\it An almost \complex structure} $\J$ is a section of fibre bundle SO$(\TT)$ with $\J^2=-$id.
Let $L_\J$ be the $(-\sqrt{-1})$-eigenspace with respect to $\J$ and $\ol L_\J$ its complex
conjugate. Then
an almost \complex structure $\J$ gives rise to the decomposition of 
the complexified $(\TT)^{\C}$ into eigenspaces 
$$
(\TT)^{\C}=L_\J\oplus \ol L_\J.
$$
An almost \complex structure $\J$ is {\it integrable} if the space $L_\J$ is involutive with respect to the Courant bracket 
$[\, ,\,]_c$. 
An integrable $\J$ is called a {\it \complex structure}. 
\subsection{$\J$-submanifolds}
Let $i_{\M} : M \to X$ be  a submanifold of $2m$ dimension.
We denote by $T^*X|_M$ the restricted bundle $i_{\M}^{-1}T^*X$ over $M$.
Let $p$ be a bundle map defined by the pull back $i_{\M}^*$ and the identity map id$_{TM}$ of $TM$,
 $$
p=\text{\rm id}_{\ss TM}\oplus i_{\M}^* : TM\oplus T^*X|_M\to 
 TM\oplus T^*M.
 $$ 
 We denote by $N^*(=N^*_{M|X})$ the conormal bundle  to $M$ in $X$.
 Then we have the short exact sequence, 
 $$
 0\arrow N^*\arrow TM\oplus T^*X|_M\overset{p}{\arrow}
 TM\oplus T^*M\arrow 0.
 $$
 We define an intersection  $L_\J(M)$ by 
 $$L_\J(M)=L_\J\cap (TM\oplus T^*X|_M)^\C=L_\J\cap (\pi^{-1}(TM))^\C$$ 
 and denote by $\ol L_\J(M)$ its complex conjugate. 
 For simplicity, we assume that $L_\J(M)$ is a subbundle of $(TM\oplus T^*X|_M)^\C$.
 Then the map $p$ is restricted to the direct sum $L_\J(M)\oplus \ol L_\J(M)$ and we have the map $q :L_\J(M)\oplus\ol L_\J(M)
 \arrow (TM\oplus T^*M)^\C$. Let $L_\J(N^*)$ denotes the intersection,
 $$
 L_\J(N^*)=L_\J\cap(N^*)^\C.
 $$
 We also assume that $L_\J(N^*)$ a subbundle of $(N^*)^\C$.
 Then $L_\J(N^*)\oplus \ol L_\J(N^*)$ is a bundle in the kernel of the map $q$ and we have the following sequence, 
 \bgn{equation}\label{eq: seq}
L_\J(N^*)\oplus \ol L_\J(N^*)\arrow  L_\J(M)\oplus\ol L_\J(M)
 \overset{q}{\arrow} (TM\oplus T^*M)^\C.
 \end{equation}
The sequence is not exact in general. 
Note that 
$$L_\J(N^*)\oplus \ol L_\J(N^*)\subset \(L_\J(M)\oplus\ol L_\J(M)\)\cap (N^*)^\C=\ker q.$$ 
 The following definition is same as those in \cite{OM}, which was given in terms of the pulback of Dirac structures \cite{Co}.
 \bgn{definition}\label{def: J-sub}
 A submanifold $M$ is a $\J$-submanifold if the sequence (\ref{eq: seq}) 
 is exact, 
 \\
 \bgn{equation}\label{eq: J-sub}
0\arrow L_\J(N^*)\oplus \ol L_\J(N^*)\arrow  L_\J(M)\oplus\ol L_\J(M)
 \overset{q}{\arrow} (TM\oplus T^*M)^\C\arrow 0.
 \end{equation}
 \smallskip
 \end{definition}
 The image $q(L_\J(M))$ is a maximally isotropic subbundle of $TM\oplus T^*M$ 
 (see \cite {Co}). Hence rank $\L_\J(M)$ =dim$_\R M$ and we have 
 \bgn{lemma}\label{lem: three equiv}
 There are three equivalent conditions: 
 \bgn{enumerate}
 \renewcommand{\labelenumi}{(\arabic{enumi})}
 \item$M$ is a $\J$-submanifold. \\
 \item The map $q :L_\J(M)\oplus\ol L_\J(M)\arrow (TM\oplus T^*M)^\C$ is surjective.\\
 \item \footnote{The condition (3) is the definition of  a \complex submanifold in \cite{OM}. 
 There is a presentation by classical tensor fields \cite{V}.}
 $q(L_\J(M))\cap q(\ol{L_\J(M)})=\{0\}$.
 \end{enumerate}
 \end{lemma}
 Since both bundles $L_\J(N^*)\oplus \ol L_\J(N^*)$ and  $L_\J(M)\oplus\ol L_\J(M)$ are $\J$-invariant, $TM\oplus T^*M$ inherits 
 the almost \complex structure $\J_{\M}$ induced in the quotient bundle.
 The almost \complex structure $\J_{\M}$ gives the decomposition into eigenspaces, 
 $$
 (TM\oplus T^*M)^\C=L_{\J,\M}\oplus\ol L_{\J, \M},
 $$
 and we have the exact sequence, 
 $$
 0\arrow L_\J(N^*)\arrow L_\J(M)\arrow L_{\J, \M}\arrow 0,
 $$
 where $L_{\J, \M}=q(L_\J(M))$.
 Then from the view point of the Dirac structure, it is shown in \cite {MO} that 
 \bgn{theorem}\cite{OM}\label{th: GCsub}
 The induced structure $\J_{\M}$ is integrable and $M$ inherits a \complex structure.
 \end{theorem}
 \bgn{proof}
For the sake of reader, we will give a proof.
We denote by $E|_M\in \Gam(M, \TT|_M)$ the restriction of a smooth section $E\in \Gam(X, \TT)$ to 
 $M$.
 Let $E_i$ be a section of $\Gam(X, \TT)$ with $E_i|_M\in \Gam (M, TM\oplus TX^*|_M)$ for $i=1,2.$
 The Courant bracket is given by 
 \bgn{equation}\label{eq: Courant}
 [E_1, E_2]_c=[u_1, u_2]+\frac12\{\L_{v_1}\theta^2-\L_{v_2}\theta^1
 -i_{v_2}d\theta_1+i_{v_1}d\theta^2\},
 \end{equation}
 where $E_i=u_i+\theta^i$ for $u_i\in TM$ and $\theta^i\in T^*X$ 
 ($i=1,2$). 
 Since the pull back $i_{\M}^*$ commutes with the exterior derivative $d$,
the Lie derivative $\L_u$ and the interior product $i_u$ for $u\in TM$, 
the Courant bracket satisfies the following, 
 \bgn{equation}\label{eq Courant}
  \bigl[\,p(E_1|_M), \,p(E_2|_M)\,\bigr]_c=p\([\,E_1, E_2]_c|_M\,\).
 \end{equation}
Since $M$ is a $\J$-submanifold, we have the exact sequence, 
$$
0\arrow L_\J(N^*)\oplus\ol L_\J(N^*)\arrow L_\J(M)\oplus \ol L_\J(M)
\overset{}{\arrow} (TM\oplus T^*M)^{\C}\arrow 0,
$$
where $(TM\oplus T^*M)^{\C}=L_{\J, \M}\oplus\ol L_{\J,\M}$. 
For every smooth section $\til E\in \Gam(M,L_{\J,M})$, there exists a smooth section $E\in \Gam(X,\TT)$ with $E|_M\in \Gam(M, L_\J(M))$ satisfying $p(E|_M)=\til E$. It follows from (\ref{eq: Courant}) that for $\til E_1,\, \til E_2\in\Gam(M,L_{\J, \M})$ we have 
$$
[\til E_1, \til E_2]_c=p\([E_1, E_2]_c|_M\).
$$
Since $[E_1, E_2]_c \,|_M\in L_\J\cap( \pi^{-1}(TM))^\C=L_\J(M)$ for $\til E_1,\, \til E_2\in \Gam(M,L_\J(M))$, we have 
$$[\til E_1, \til E_2]_c\in p(L_\J(M))=L_{\J,\M}.$$
Hence $\J_{\M}$ is integrable.
\end{proof}
Note that theorem \ref{th: GCsub} holds for sheaves $L_\J(M)$ and $L_\J(N^*)$ with a similar proof.
\bgn{example}(symplectic submanifolds)
A symplectic manifold with a symplectic structure $\ome$ gives the 
\complex structure $\J_\ome$. Then
a symplectic submanifold of a symplectic manifold is a $\J_\ome$-submanifold. 
In this case we have 
$$
L_{\J_\ome}(M)\oplus \ol L_{\J_\ome}(M) \cong (TM\oplus T^*M)^\C,
$$
with $L_{\J_\ome}(N^*)=\{0\}$.
\end{example} 
\bgn{example}(complex submanifolds)
A complex manifold with a complex structure $J$ admits the \complex structure $\J_J$. Then
a complex submanifold of a complex manifold is a $\J_J$-submanifold with the exact sequence 
\bgn{equation*}
0\arrow L_{\J_J}(N^*)\oplus \ol L_{\J_J}(N^*)\arrow  L_{\J_J}(M)\oplus\ol L_{\J_J}(M)
 \overset{q}{\arrow} (TM\oplus T^*M)^\C\arrow 0,
 \end{equation*}
 where $L_{\J_J}(N^*)=(N^*)^{1,0}$ and $L_{\J_J}(M)=TM^{0,1}\oplus (T^*X)^{1,0}|_M$.
\end{example}
\bgn{example}
Let $X_1\times X_2$ be the product of \complex manifolds $(X_1, \J_1)$ and $(X_2,\J_2)$ with the \complex structure $\J_1\times \J_2$.
Let $M_i$ be a $\J_i$-submanifold of $(X_i, \J_i)$ for $i=1,2$. 
Then the product $M_1\times M_2$ is a $(\J_1\times \J_2)$-submanifold 
of $X_1\times X_2$.
\end{example}
\bgn{example}\label{ex: invariant conormal}($\J$-invariant conormal bundles)
Let $(X, \J)$ be a \complex manifold and $i_{\M}\: M\to X$ a submanifold of $X$ whose conormal bundle $N^*$ is $\J$-invariant, 
$$
\J(N^*)=N^*.
$$
Then $M$ is a $\J$-submanifold. \\
\bgn{proof}{\it of example \ref{ex: invariant conormal}}
The bundle $TM\oplus T^*X|_M$ is also defined in terms of $N^*$, 
$$
TM\oplus T^*X|_M=\{\, E\in TX\oplus T^*X|_M\,\,
\lan E, \theta\ran =0, \forall \theta\in N^*\,\}.
$$
Since $\J$ is a section of SO$(\TT)$, we have 
$$
\lan \J E, \theta\ran =-\lan E, \J \theta\ran =0,
$$
for $E\in TM\oplus T^*X|_M$ and $\theta\in N^*$. Thus the bundle 
$TM\oplus T^*X|_M$ is $\J$-invariant.
Hence we have the exact sequence by 
\[\xymatrix{
0\ar[r]& (N^*)^\C\ar[r] \ar@{=}[d]&(TM\oplus T^*X|_M)^\C\ar[r] \ar@{=}[d]&(TM\oplus T^*M)^\C\ar[r]\ar@{=}[d]&0\\
0\ar[r]& L_\J(N^*)\oplus \ol L_\J(N^*)\ar[r] &L_\J(M)\oplus \ol L_\J(M)\ar[r] &(TM\oplus T^*M)^\C\ar[r]&0
}\]
\end{proof}
\end{example}
\bgn{example}\label{ex: b-fields}($b$-fields)
Let $M$ be a $\J$-submanifold of a \complex manifold $(X,\J)$. 
For a real $d$-closed $2$-form $b$, the exponential $e^b$ acts on $\J$ by the adjoint action to obtain a \complex structure  $\J_b=\Ad_{e^b}\J$. 
Then $M$ is also a $\J_b$-submanifold.
\end{example}
\bgn{proof}{\it of example \ref{ex: b-fields}}
The bundles $L_{\J_b}(M)$ and $L_{\J_b}(N^*)$ are  respectively
given as the images by the adjoint action $\Ad_{e^b}$,
$$
L_{\J_b}(M)=\Ad_{e^b}\( L_\J(M)\),\qquad 
L_{\J_b}(N^*)=\Ad_{e^b}\(L_\J(N^*)\),
$$
and the adjoint action $\Ad_{e^b}$ preserves the bundles $TM\oplus T^*X|_M$ and $N^*$. Hence we have the exact sequence by 
\[\xymatrix{
0\ar[r]& L_\J(N^*)\oplus \ol L_\J(N^*)\ar[r] \ar@{=}[d]&(L_\J(M)\oplus \ol L_\J(M))\ar[r] \ar[d]^{\Ad_{e^b}}&(TM\oplus T^*M)^\C\ar[r]\ar[d]^{\Ad_{e^{i^*b}}}&0 \\
0\ar[r]& L_{\J_b}(N^*)\oplus \ol L_{\J_b}(N^*)\ar[r] &(L_{\J_b}(M)\oplus \ol L_{\J_b}(M))\ar[r] &(TM\oplus T^*M)^\C\ar[r]&0 ,
}\]
where $\Ad_{e^{i^*b}}$ denotes the adjoint action by the  exponential of the pull back $i^*b$.
\end{proof}
 \subsection{Generalized metrics and generalized K\"ahler structures}
{\it A generalized metric} $\hat G$ is a section of SO$(\TT)$ with $\h G^2=$id satisfying the condition: a bilinear form $G$ defined by 
$
G(E_1, E_2):= \lan \h G E_1, E_2\ran $
is  a positive-definite metric on $\TT$, where 
$E_i=v_i+\eta_i$ for $v_i\in TX$ and $\eta_i\in T^*X$ ($i=1,2$). 
A \metric gives the decomposition of $\TT$ into eigenspaces 
$$
\TT=C^+\oplus C^-,
$$
where $C^+$ and $C^-$ denotes the $(+1)$-eigenspace and the $(-1)$-eigenspace respectively. 
Then there are a Riemannian metric $g$ and a $2$-form $b$ such that 
$C^+$ and $C^-$ are respectively written as 
\bgn{align}\label{eq: Cpm}
C^+ =&\{ \, v + g( v,\, )+b(v,\,)\, |\, v\in T\,\}\\
C^-=&\{ \, v - g( v,\, )+b(v,\,)\, |\, v\in T\,\},\nonumber
\end{align}
where $g(v, \,)$ and $b(v, \,)$ denote the $1$-forms given by the interior product by $v\in TX$ respectively. Hence the restriction of the projection $\pi$ to $C^+$ and $C^-$
gives isomorphisms respectively,
\bgn{equation}\label{C pi}
\pi|_{C^+} : C^+ \cong TX, \quad \pi|_{C^-} : C^-\cong TX
\end{equation}
\\
{\it A \K\"ahler structure } is a pair $(\J_0, \J_1)$ consisting of two commuting \complex structures with the \metric $G=-\J_0\J_1=-\J_1\J_0$.
Then $(\TT)^\C$ is simultaneously decomposed into
four eigenspaces by $\J_0$ and $\J_1$, 
\bgn{align}
(\TT)^\C=&(L_{\J_0}\cap L_{\J_1})\oplus(L_{\J_0}\cap\ol L_{\J_1})\\
\oplus&
(\ol L_{\J_0}\cap L_{\J_1})\oplus(\ol L_{\J_0}\cap\ol L_{\J_1}).
\end{align}
Then eigenspaces $C^+$ and $C^-$ of the \metric $G=-\J_0\J_1=-\J_1\J_0$ are respectively given by 
\bgn{align}
(C^+)^\C=&(L_{\J_0}\cap L_{\J_1})\oplus(\ol L_{\J_0}{\cap}\ol L_{\J_1}),\\
(C^-)^\C=&(L_{\J_0}{\cap}\ol L_{\J_1})\oplus(\ol L_{\J_0}{\cap}L_{\J_1})
\end{align}
\subsection{Generalized K\"ahler submanifolds}
\label{subsec: GKsub}
As in example \ref{ex: invariant conormal}, if $M$ has a $\J$-invariant conormal bundle, 
then $M$ admits the induced \complex structure $\J_{\M}$. Then we shall show the following in this section 
\ref{subsec: GKsub}
\bgn{theorem}\label{th: GKsub}
 Let $(\J_0,\J_1)$ be a \K\"ahler structure on $X$. 
If a submanifold $M$ of $X$ admits a $\J_0$-invariant conormal bundle, then $M$ is also a $\J_1$-submanifold and 
$M$ inherits the induced \K\"ahler structure $(\J_{0,\M}, \J_{1,\M})$.
\end{theorem}
\bgn{remark}
 Gualtieri pointed out that theorem \ref{th: GKsub}
can be proved by a modified version of the generalized K\"ahler reduction \cite{B.G.C}.
\end{remark}
At first we define a \K\"ahler submanifold which is due to Barton and Sti{\' e}non \cite{BS}.
\bgn{definition}
Let $(X, \J_0, \J_1)$ be a \K\"ahler manifold and 
$M$  a submanifold of $X$. 
A submanifold $M$ is a \K\"ahler submanifold if $M$ is a $\J_0$-submanifold and $M$ is also a $\J_1$-submanifold.
\end{definition}
Then it is shown that 
a \K\"ahler submanifold $M$
 inherits a \K\"ahler structure $(\J_{0,\M}, \J_{1,\M})$. \\ \indent
Let $(X, \J_0, \J_1) $ be a $2n$-dimensional \K\"ahler manifold  with the \metric $G$ and
$M$ a submanifold of dimension $2m$.
As in section 1.2, we have subbundles $C^+$ and $C^-$. 
We define $C^+(M)$ and $C^-(M)$ respectively by the intersections
\bgn{equation}
C^+(M) =C^+\cap \pi^{-1}(TM), \quad C^-(M)= C^-\cap\pi^{-1}(TM).
\end{equation} 
Then from (\ref{C pi}), we see that $C^+(M)$ and $C^-(M)$ are bundles with
$\rank\, C^+(M)=\rank\, C^-(M)=\dim M=2m$.
Let $p$ be the bundle map in section 1.2,
$$
p : TM\oplus T^*X|_M\arrow TM\oplus T^*M
$$
We denote by $\gam$ the bundle map given by the restriction of $p$ to the subbundle 
$(C^+(M)\oplus C^-(M))$,
\bgn{lemma}\label{lem: gamiso}
The map $\gam :  (C^+(M)\oplus C^-(M))\arrow  (TM\oplus T^*M)$ is an isomorphism.
\end{lemma}
\bgn{proof}
The kernel of the map $\gam$ is the intersection $N^*\cap  (C^+(M)\oplus C^-(M))$.
From (\ref{eq: Cpm}), $C^+(M)$ and $C^-(M)$ are respectively written as
\bgn{align}
C^+(M)=&\{\,u_1+g(u_1, \,)+b(u_1, \,)\, |\, u_1\in TM\,\},\\
C^-(M)=&\{\,u_2-g(u_2, \,)+b(u_2, \,)\, |\, u_2\in TM\,\}.
\end{align}
Since $g$ is positive-definite, it follows that $N^*\cap (C^+(M)\oplus C^-(M))=\{0\}.$
Hence $\gam$ is injective. Since $\rank\,(C^+(M)\oplus C^-(M))=\rank\,(TM\oplus T^*M)=4m$, it follows that $\gam$ is an isomorphism.
\end{proof}
\bgn{lemma}\label{lem: gamH(M)}
If $M$ admits a $\J_0$-invariant conormal bundle, then $C^+(M)$ and $C^-(M)$ are respectively invariant under both action of $\J_0$ and $\J_1$.
\end{lemma}
\bgn{proof}
If $M$ admits a $\J_0$-invariant conormal bundle
as in example \ref{ex: invariant conormal}, $\pi^{-1}(TM)=TM\oplus T^*X|_M$ is also 
$\J_0$-invariant. Since $C^+$ and $C^-$ are respectively eigenspaces of $G$ and $\J_0$ commutes with $G$, $C^+(M)$ and $ C^-(M)$ are respectively invariant under both action $\J_0$ and $G$. It follows from $\J_1=G\J_0$ that $\J_1$ is preserving $C^+(M)$ and $ C^-(M)$. 
\end{proof}\\
\bgn{proof}{\it of theorem \ref{th: GKsub}}
We shall show the sequence is exact,
$$
L_{\J_1}(N^*)\oplus\ol L_{\J_1}(N^*)\arrow 
L_{\J_1}(M)\oplus \ol L_{\J_1}(M) \overset{q_1}{\arrow}(TM\oplus T^*M)^\C.
$$
Since $M$ is a $\J_0$-submanifold, it follows from theorem \ref{th: GCsub} that we have the exact sequence with respect to $\J_0$,
$$
0\arrow L_{\J_0}(N^*)\oplus \ol L_{\J_0}(N^*)
\arrow L_{\J_0}(M)\oplus \ol L_{\J_0}(M)\overset{q_0}{\arrow} (TM\oplus T^*M)^\C 
\arrow 0,
$$
and we have the induced \complex structure $\J_{o,\M}$.
From lemma \ref{lem: gamH(M)}, $(C^+(M)\oplus C^-(M))^\C$ is a subbundle of both $L_{\J_0}(M)\oplus \ol L_{\J_0}(M)$ and $L_{\J_1}(M)\oplus \ol L_{\J_1}(M)$.
It follows from lemma \ref{lem: gamiso} that we have the following commutative diagram, 
\[\xymatrix{
0\ar[r] &L_{\J_0}(N^*)\oplus\ol L_{\J_0}(N^*)\ar[r]&L_{\J_0}(M)\oplus
\ol L_{\J_0}(M)\ar[r]^{q_0}&(TM\oplus T^*M)^\C\ar[r]&0\\
&&(C^+(M)\oplus C^-(M))^\C\ar[u]\ar[d]\ar[r]^{\cong}_\gam&(TM\oplus T^*M)^\C\ar@{=}[u]\ar@{=}[d]&\\
 &L_{\J_1}(N^*)\oplus\ol L_{\J_1}(N^*)\ar[r]&L_{\J_1}(M)\oplus
\ol L_{\J_1}(M)\ar[r]_{q_1}&(TM\oplus T^*M)^\C&
}
\]
Hence the map $q_1 :L_{\J_1}(M)\oplus
\ol L_{\J_1}(M)\to(TM\oplus T^*M)^\C$ is surjective. 
Hence it follows from lemma \ref{lem: three equiv} (2) that $M$ is a $\J_1$-submanifold and a \K\"ahler submanifold.
\end{proof} \\
Our theorem 1.9 can be generalized. For instance,
as in our proof, if $C^+(M)\oplus C^-(M)$ is invariant under the action of $\J_0$, then 
$M$ is a \K\"ahler submanifold.

\subsection{Generalized K\"ahler manifolds with one pure spinor}
{\it A pure spinor} of $X$ is  a complex differential form $\psi$ with 
$\dim_\C\ker\psi=2\dim_\C X,$ where $\ker\psi=\{\, E\in (\TT)\otimes\C\, |\, 
E\cdot \psi=0\, \}$.
A pure spinor is {\it non-degenerate} if we have a decomposition,
$$
(\TT)^\C=\ker\psi\oplus \ol{\ker\psi}.
$$
Thus a non-degenerate, pure spinor induces the almost \complex structure $\J_\psi$ such that $\ker\psi$ is the $(-\sqrt{-1})$ eigenspace $L_\psi$ of $\J_{\psi}$. 
If a non-degenerate, complex pure spinor $\psi$ is $d$-closed, then the induced $\J_\psi$ is integrable. 
\bgn{definition}
A pair $(\J_0, \psi)$ consisting of a \complex structure and a $d$-closed, non-degenerate, complex pure spinor is a \K\"ahler structure with one pure spinor if the induced pair $(\J_0, \J_\psi)$ 
is a \K\"ahler structure.
\end{definition}
For a point $x\in M$, the pure spinor $\psi$ is written as 
\bgn{equation}\label{eq: pure spinor}
\psi_x=\psi_{l,x}e^{b+\sqrt{-1}\ome},
\end{equation}
where $\psi_{l,x}$ is a complex $l$-form which is given by 
$\psi_{l, x}=\theta^1\w\cdots\w\theta^l$ in terms of  $1$-forms $\{\theta^i\}_{i=1}^l$
 and $\ome$ and $b$ are real $2$-forms.
 The degree of $\psi_{l,x}$ is called Type of a pure spinor $\psi$ at $x$. 
 If Type $\psi_x=0$, it follows that 
 the pullback $i_{\M}^*\psi_x$ of $\psi$ to a submanifold $M$ does not vanish.
 In general the pullback $i_{\M}^*\psi$ may vanish which is not a pure spinor on $M$.
 However, we have 
\bgn{theorem}\label{th: GKPsub}
Let $(X, \J_0, \psi)$ be a \K\"ahler manifold with one pure spinor.
Let $M$ be a submanifold with invariant conormal bundle with respect to $\J_0$.
Then the pull back $i_{\M}^*\psi$ is 
a $d$-closed, non-degenerate, pure spinor on $M$ and the induced pair $(\J_{0,\M}, i_{\M}^*\psi)$ is 
a \K\"ahler structure with one pure spinor on $M$.
\end{theorem}
We shall show the following lemma for proof of theorem \ref{th: GKPsub}.
\bgn{lemma}\label{lem: GKPsub} 
Let $(X,  \J_0, \psi)$ be a \K\"ahler manifold with one pure spinor and 
$M$ a submanifold with invariant conormal bundle with respect to $\J_0$. Then the pull back $i_{\M}^* \psi $ does not vanish.
\end{lemma}
\bgn{proof}{\it of lemma \ref{lem: GKPsub}}
 In the case $l=0$, then $i_{\M}^*\psi_x=i_{\M}^*e^{b+\sqrt{-1}\ome}\neq 0$.
 Thus it suffices to consider the case $l>1$. 
 From (\ref{eq: pure spinor}), if $i_{\M}^*\psi_x=0$, then we have $i_{\M}^*\psi_{l,x}=0$.
 Thus $i_{\M}^*\psi_{l,x}$ is generated by $N^*$ and  at least one element of $\{\theta^i\}_{i=1}^l$ belongs to $N^*$.  We can assume that $\theta^i\neq0\in N^*$. 
It follows from $\theta^i\cdot \psi_x=0$ that we see $\theta^i\in L_\psi$.
Then we have 
\bgn{align}
G(\theta^i, \ol\theta^i)=&\lan G\theta^i, \ol\theta^i\ran=
\lan -\J_0\J_\psi \theta^i, \ol\theta^i\ran\\
=&\sqrt{-1}\lan \J_0\theta^i, \ol\theta^i\ran.
\end{align}
Since $\J_0 N^*=N^*$, we see that $\J_0\theta^i$ is a $1$-form and 
$\lan \J\theta^i, \ol\theta^i\ran=0$. 
Hence $G(\theta^i, \ol\theta^i)=0$. Since $G$ is positive-definite, it implies $\theta^i=0$,
which is a contradiction. Hence we conclude that $i_{\M}^*\psi_x\neq0$ for all $x\in M$.
\end{proof}

\bgn{proof}{ \it of theorem \ref{th: GKPsub}}
Let $\J_\psi$ be the induced \complex structure by $\psi$ with the $(-\sqrt{-1})$-eigenspace 
$L_\psi$. As in section 1.4, $C^+(M)\oplus C^-(M)$ is $\J_\psi$-invariant and 
under the isomorphism $q : C^+(M)\oplus C^-(M)\cong TM\oplus T^*M $ we have the decomposition, 
$$
(TM\oplus T^*M )^\C = q( L_\psi(M))\oplus q(\ol L_\psi(M)).
$$
For $E= u +\eta\in L_\psi^M$, we see that 
\bgn{align}
q(E)\cdot i_{\M}^*\psi =&(u+i_{\M}^*\eta)\cdot i_{\M}^*\psi\\
 =&i_{\M}^*( u+\eta)\cdot\psi= i_{\M}^*( E\cdot \psi)=0.
\end{align}
It implies that $q(L_\psi(M))\subset \ker i_{\M}^*\psi$. 
Since $\dim q(L_\psi(M))=2m$ and $q(L_\psi(M))\cap q(\ol L_\psi(M))=\{0\}$, 
it follows from lemma \ref{lem: GKPsub}  that $q(L_\psi(M))=\ker i_{\M}^*\psi$ is maximally isotropic.
Thus $i_{\M}^*\psi$ is a non-degenerate, pure spinor on $M$ with 
the induced structure $\J_{\M,\psi}$.
Since the pull back $i_{\M}^*\psi$ is $d$-closed, the pair 
$(\J_{\M,0}, \J_{\M,\psi})$ is a \K\"ahler structure. 
Hence the pair $(\J_{\M,0}, i_{\M}^*\psi)$ is a \K\"ahler structure with one pure spinor.
\end{proof}
\subsection{Poisson submanifolds}
\bgn{definition}
Let $X$ be a complex manifold with a holomorphic $2$-vector $\b$. 
If the Schouten bracket vanishes, i.e., $[\b,\b]_{Sch}=0$, we call $\b$  a (holomorphic) Poisson structure on $X$ and the Poisson bracket is defined by 
$$
\{ f, g\} =\b(df\w dg).
$$
\end{definition}
\bgn{definition}
Let $X$ be a complex manifold with a Poisson structure $\b$ and $M$ a complex submanifold with the defining ideal sheaf $I_M$. 
A submanifold $M$ is a Poisson submanifold if we have $\{ f, g\} \in I_M$ for all 
$f\in I_M$ and $g\in {\cal O}_X$.
\end{definition}
A Poisson submanifold admits the induced Poisson structure.\\
Let $\J_J$ be the generalized complex structure defined by the usual complex structure $J$. 
By using a Poisson structure $\b$, we obtain a family of generalized complex structures $\J_{\b,t}$ parameterized by $t$
$$
\J_{\b t} =e^{a t}\circ \J_J\circ e^{-a t},
$$
where $a=\b+\ol\b$.
The structure $\J_{\b t}$ is written in the form of a matrix, 
\bgn{equation}\label{eq: Jb}
\J_{\b t}=\bgn{pmatrix}
J&-Ja-aJ^*\\
0&-J^*
\end{pmatrix}.
\end{equation}
Then we have 
\bgn{lemma}\label{poisson sub}
Every Poisson submanifold $M$ admits a $\J_{\b t}$-invariant conormal bundle for all $t$.
\end{lemma}
\bgn{proof}
Since $\b( df, dg) \in I_M$ for $f\in I_M$ and $(N^*)^{1,0}$ is generated by the set $\{ df\, |\,f\in I_M\}$,
we have the restriction $\b(df, \, )|_M=0$.
It follows from (\ref{eq: Jb}) that 
 $\J_\b (N^*)=N^*$.
\end{proof}\\
In \cite{Go3} we obtain a stability theorem of \K\"ahler structure with one pure spinor. 
It implies that \K\"ahler structure with one pure spinor is stable under small deformations
of \complex structures. By applying the stability theorem to small deformations of 
\complex structures $\{ \J_{\b t}\}$ starting from $\J_J$, we have 
deformations of \K\"ahler structures with one pure spinor $\{\J_{\b t}, \psi_t\}$. 
The type of $\J_{\b t}$ is given by 
$$
\text{\rm Type}\J_{\b t} = n -2 \text{ \rm rank}\, \b.
$$
It implies that if $\b\neq 0$, then deformations of \K\"ahler structures with one pure spinor $\{\J_{\b t}, \psi_t\}$ can not be obtained from ordinary K\"ahler structures by the action of $b$-fields. 
Hence we have,
\bgn{theorem}\label{th: stability}
Let $X$ be a K\"ahler manifold with non-trivial Poisson structure $\b$. 
Then there exists an analytic family of non-trivial \K\"ahler structure with one pure spinor 
$\{\J_{\b,t}, \psi_t\}$.
\end{theorem}
Hence it follows from \ref{poisson sub} and \ref{th: GKPsub} that 
\bgn{theorem}\label{th:Poisson sub}
Let $M$ be a Poisson submanifold of a Poisson manifold $X$ with a K\"ahler structure $\ome$. Then $M$ is a \K\"ahler manifold with the induced structure 
$(\J_{\b t,\M}, \psi_{t,\M})$.
\end{theorem}
\subsection{Examples of \K\"ahler submanifolds arising as
Poisson submanifolds}
Let $X$ be a compact K\"ahler manifold on which an $l$ dimensional commutative complex group $G$ act holomorphically. 
The Lie algebra $\frak g$ of $G$ generates holomorphic vector fields $\{ V_i\}_{i=1}^l$ on $X$. 
Since $[V_i, V_j]=0$, it follows that a linear combination of $2$-vectors $V_i\w V_j$'s gives a holomorphic Poisson structure $\b$, 
\bgn{equation}\label{eq: Poisson by commutative action}
\b=\sum_{i,j} \lam_{i,j} V_i\w V_j,
\end{equation}
where $\lam_{i,j}$ is a constant.
Note that $[V_i, V_j]=0$ implies $[\b, \b]_{Sch}=0$.
If $\b\neq0 $, from the stability theorem we have  deformations of \K\"ahler structures 
$\{\J_{\b t}, \psi_t\}$. 
\bgn{lemma}
Let $M$ be a complex submanifold with defining ideal $I_M$. 
If the ideal $I_M$ is invariant under the action of $G$, then $M$ is a Poisson submanifold of 
$(X, \b)$. 
\end{lemma}
\bgn{proof}
Since we have $V_if\in I_M$ for $f\in I_M$ and $i=1,\cdots ,l$.
Hence
$\b( df,\, )\in I_M\otimes T^{1,0}X$ for $f\in I_M$. 
It implies that $M$ is a Poisson submanifold.
\end{proof}
\bgn{example}(toric submanifolds)
Let $X$ be a compact toric manifold of $n$ dimension. Then there is the action of 
$n$ dimensional complex torus $G$ on $X$. Then a toric submanifold $M$ which is invariant under the action of $G$ is a Poisson submanifold with respect to $\b$ as in (\ref{eq: Poisson by commutative action}).
\end{example}
\bgn{example}
Let $\C P^4$ be the complex projective space with the homogeneous coordinates $[z_0, z_1, z_2, z_3, z_4]$ on which the commutative group 
$\C^\times \times \C^\times$ act by a homomorphism $\rho : \C^\times \times \C^\times\to $GL$(5, \C)$, 
$$
\rho(\lam_1, \lam_2) =\text{\rm diag}(1, \lam_1, \lam_1, \lam_2, \lam_2).
$$
Then we have a Poisson structure $\b=V_1\w V_2$ as in (\ref{eq: Poisson by commutative action}).
We take a following quadratic function $F$ of 
$\C P^4$
$$
F(z)=\sum_{\overset{ i=1,2}{\ss j=3,4}} a_{i j} z_i z_j,
$$
where $a_{ij}$ are constants.
The hypersurface $M$ defined as the zero of $F$ becomes a smooth manifold of complex $3$ dimension for suitably choosed constants  $a_{ij}$.
Since $\rho ^*(\lam_1, \lam_2) F(z) =\lam_1\lam_2F(z)$, the hypersurface $M$ is a Poisson submanifold in $\C P^4$ 
which admits the deformations \K\"ahler structure with one pure spinor 
$(\J_{\b t}, \psi_t)$ from theorem \ref{th:Poisson sub}. Since the Type of $\J_{\b t}$ 
is $2$ at generic points of $X$ and the type of the induced $\J_{\b t, \M}$ is 1 at generic points of $M$, the \K\"ahler structure 
$(\J_{\b t}, \psi_t)$ and the induced \K\"ahler structure 
$(\J_{\b t \M}, \psi_t)$ are not obtained from 
K\"ahle structures by the action of $b$-fields.
\end{example}
\bgn{example}
Let $\C P^3$ be the complex projective space with 
the homogeneous coordinates $[z_0, z_1, z_2,z_3]$. 
On the open set $\{z_0\neq 0\}$, we have the inhomogeneous coordinates 
$\{\zeta_1, \z_2, \zeta_3\}$ given by 
$\displaystyle{\z_i=\frac{z_i}{z_0}}$. Let $f=f(\zeta_1, \z_2, \zeta_3)$ be a polynomial of degree $d\leq 3$ and we assume that  $df \neq 0.$
Then we define a $2$-vector field $\b_f$ by 
$$
\b=f_1\frac{\pa}{\pa \z_2}\w\frac{\pa}{\pa \z_3}+
f_2\frac{\pa}{\pa \z_3}\w\frac{\pa}{\pa \z_1}+ f_3\frac{\pa}{\pa z_1}\w\frac{\pa}{\pa \z_2},
 $$where $\displaystyle{f_i=\frac{\pa}{\pa z_i}f}$. Then it turns out that 
 $[\b ,\b]_{\sch}=0$. Thus $\b_f$ is a Poisson structure, which is called the exact quadratic Poisson structure \cite{L.P},\cite{Pol}.
 We also see that $\b(df)=0$. 
 Thus the zero of $f$ is a Poisson submanifold with respect to 
 the Poisson structure $\b_f$ on $\C^3$. 
 Let $F=F(z_0, \cdots, z_3)$ be the homogeneous polynomial defined by 
 $$F=z_0^df\(\frac{z_1}{z_0}, \frac{z_2}{z_0}, \frac{z_3}{z_0}\).$$
Since  each $f_i$ is a quadratic polynomial, 
$\b_f$ can be extended as a holomorphic Poisson
structure $\b_F$ on $\C P^3$. Then a complex surface $M$ given by the zero of $F$ is a Poisson submanifold.
 \bgn{theorem}
 Let $M$ be a complex smooth hypersurface of the projective space $\C P^3$
 defined by a homogeneous polynomial $F$ of degree $d\leq 3$ .
 Then $M$ is a non-trivial \K\"ahler manifold arising as  Poisson submanifold of $\C P^3$ with respect to the Poisson structure $\b_F$. 
 \end{theorem}
 \bgn{proof}
 It suffices to show that the induced Poisson structure $\b_{\M}$ is non-trivial. 
 On $\{ z_0\neq 0\}$, $\b_{\M}$ is the induce structure from $\b_f$.
 Since $M$ is smooth, we can assume that there exists an open set defined by
 $\{f_3\neq 0\}$
 with coordinates $(\eta_1, \eta_2, \eta_3)$, 
 \bgn{align*}
 &\eta_1=\z_1\\
& \eta_2=\z_2\\
 &\eta_3=f(\z_1,\z_2,\z_3).
 \end{align*}
 Then $\b_f$ is written as 
 $$
 \b_f=f_3\frac{\pa}{\pa \eta_1}\w\frac{\pa}{\pa \eta_2}.
 $$
 Since $M$ is defined by $\eta_3=0$, it follows that $\b_{\M}$ is non-trivial.
 Hence the type of the induced \complex structure $\J_{\M}$ is $0$ on 
 the complement of the zero of $\b_{\M}$. Hence the induced 
 \K\"ahler structure on $M$  is not obtained from ordinary K\"ahler structures by $b$-field action.
 \end{proof}
\end{example}
}


\large{   
\section{Deformations of \K\"ahler structures via Poisson structures} 
Let $X$ be a compact K\"ahler manifold with a K\"ahler form $\ome$. 
Then we have the \K\"ahler structure with one pure spinor $(\J_J, \psi)$, 
where $\J_J$ denotes the \complex structure induced from the complex structure $J$ on $X$ and 
$\psi$ is the pure spinor defined by 
$$
\psi=e^{\sqrt{-1}\ome}.
$$
We assume that there exists a Poisson structure $\b$ on $X$. Then we have deformations of 
\complex structures $\{\J_{\b t}\}_{t\in \trian}$ as in section 1. 
Applying the stability theorem in \cite{Go3} to deformations of 
\complex structures $\{\J_{\b t}\}_{t\in \trian}$, we obtain
\\\\
{\bf Theorem \ref{th: stability}}
Let $X$ be a K\"ahler manifold with non-trivial Poisson structure $\b$. 
Then there exists an analytic family of non-trivial \K\"ahler structure with one pure spinor 
$\{\J_{\b t}, \psi_t\}$.
\\\\
In the case of deformations starting from ordinary K\"ahler manifolds, proof of the stability theorem becomes simple which is based on the ordinary Hodge decomposition and Lefschetz decomposition.
Note that in general case, we used the generalized Hodge decomposition.
We shall give an exposition of our proof in the special cases.
We use the same notation as in \cite{Go3}.
\\
Let CL be the real Clifford algebra of $\TT$ with respect to $\lan\,,\,\ran$.
Then $\CL$ acts on differential forms by the spin representation.
The Clifford group $G_{cl}$ is defined in terms of the twisted adjoint $\wtil\Ad_g$,
$$
G_{cl}:=\{\, g\in \CL^\times\, |\, \wtil{\Ad}_g(\TT)\subset\TT\, \}.
$$
Let CL$^2$ be the Lie algebra which consists of elements of 
the Clifford algebra of degree less than or equal to $2$.
It turns out that CL$^2$ is the Lie algebra of the Clifford group  $G_{cl}$.
The set of almost \complex structures forms an orbit of the adjoint action of the Clifford group and 
the set of almost \K\"ahler structures with one pure spinor is also an orbit of the diagonal action of 
the Clifford group. Thus it follows that small deformations of 
almost \K\"ahler structures with one pure spinor are given by the action of exponential of CL$^2$ on 
$(\J_J, \psi)$, 
$$
\(\Ad_{e^{z(t)}} \J_J,\,\, e^{z(t)}\cdot \psi \),
$$
where $z(t) \in $CL$^2[[t]]$.
Let $\{\J_{\b t}\}$ be deformations of \complex structures by a Poisson structure $\b$ as before.
The deformations $\{\J_{\b t}\}$ are given by the adjoint action of real $2$-vector $a=\b+\ol \b$, 
$$
\J_{\b t}=\Ad_{e^{at}}\J_J.
$$
Let $b(t)$ be an analytic family of $b(t)$ of CL$^2[[t]]$, 
$$
b(t)=b_1 t+b_2\frac{t^2}{2!}+\cdots=\sum_{i=1}^\infty b_i\frac{t^i}{i!}.
$$
We denote by $\w^{n,0}$ the canonical line bundle of $(X, J)$. 
We assume that there exists a family $\{b(t)\}$ with the followings conditions 
(2.1) and (2,2),
\bgn{align}
&b_i\cdot \w^{n,0}\subset\w^{n,0}\label{condition 2.1}\\
&d (e^{at}\, e^{b(t)}\cdot\psi) =0.\label{condition 2.2}
\end{align}

Then it follows from the Campbel-Hausdorff formula that there is the $z(t)\in \CL^2[[t]]$  
with 
$$
e^{z(t)}=e^{at}e^{b(t)}.
$$
From (\ref{condition 2.1}), we see that the action by $b(t)$ is preserving $\J_J$,
$$
\Ad_{e^b(t)}\J_J =\J_J.
$$
Thus we have 
\bgn{align}
\Ad_{z(t)}\J_J=&\Ad_{at}\circ \Ad_{b(t)}\J_J\\
=&\Ad_{at}\J_J = \J_{\b t}
\end{align}
From (\ref{condition 2.2}), the non-degenerate pure spinor $\psi_t=e^{z(t)}\cdot\psi$ is $d$-closed. 
Hence the pair  $(\Ad_{z(t)}\J_J,\, e^{z(t)}\cdot\psi) =(\J_{\b t}, \psi_t)$ 
is a \K\"ahler structure with one pure spinor.
We shall construct $b(t)$ which satisfies the (\ref{condition 2.1}) and (\ref{condition 2.2}).
Let $\CL^{[i]}$ be the subspace of the $\CL$ of degree $i$. 
We define $\CL^i$ for $i=0,\cdots ,3 $ by 
\bgn{align}
&\CL^0=C^\infty(X), \quad \CL^1=\TT,\\
&\CL^2=\CL^0\oplus\CL^{[2]},\quad \CL^3=\CL^1\oplus\CL^{[3]}.
\end{align}
Then we define bundles $\wtil \ker^1$ and $\wtil\ker^2$ respectively by 
\bgn{align}
&\wtil \ker^1=
\{ \, b \in \CL^{2}\, |\, b\cdot \w^{n,0}\subset\CL^{0}\cdot\w^{n,0}\, \},\\
&\wtil \ker^2=\{ \, b \in \CL^{3}\, |\, b\cdot \w^{n,0}\subset\CL^{1}\cdot\w^{n,0}\, \},
\end{align}
 where $\CL^i\cdot\w^{n,0}$ denotes the image by the action of $\CL^i$ on the canonical line bundle $\w^{n,0}$.
Then a section $b\in \ker^i$ ($i=1,2$) acts on $\psi=e^{\sqrt{-1}\ome}$ 
by the spin representation and we obtain bundles 
$\wtil K^1$ and $\wtil K^2$, 
$$
\wtil K^i=\{\, b\cdot\psi \, |\, b\in \wtil\ker^i\, \}.
$$
The bundle $\wtil K^1$ is the direct sum of $U^{0,-n}$ and $U^{0,-n+2}$, 
$$
\wtil K^1=U^{0,-n}\oplus U^{0,-n+2},
$$
where $U^{0,-n}=\CL^0\cdot \psi=\{ \, f\psi \, |\, f\in C^\infty(X)\, \}$
and $U^{0,-n+2}$ is given by the contraction $\w_\ome$ by the K\"ahler form $\ome$,
\bgn{align}
&U^{0, -n+2}=\{\, h\psi+p\w\psi \, |\, h \in C^\infty(M), \,p\in \w^{1,1},\,\,
\w_{\ome} p +2h =0\, \},
\end{align}
where $\w^{1,1}$ denotes forms of type $(1,1)$ with respect to the complex structure $J$.
We define $K^1$ to be $U^{0, -n+2}$ and write $\wtil K^2$ as $K^2$.
Then $K^2$ is written as 
$$
K^2=\wtil K^2=\{ \, \eta\w\psi\, |\, \eta\in \w^1\oplus \w^{2,1}\oplus\w^{1,2}\, \}.
$$
Then we have a differential  complex $\{ K^i, d\}$ by the exterior derivative $d$, 
\[\xymatrix{
0\ar[r] &K^1\ar[r]^d &K^2\ar[r]^d &\cdots. 
}
\]
It turns out that the complex $(K^i, d)$ is elliptic since we have the following elliptic complex, 
\[\xymatrix{
0\ar[r] &P^{1,1}\ar[r]^d &\w^{2,1}\oplus\w^{1,2}\ar[r]^d &\cdots, 
}
\]
where $P^{1,1}$ denotes the primitive $(1,1)$-forms on the K\"ahler manifold $X$ 
({\it cf}, proposition 4.7 in \cite{Go1}).
The complex $(K^*, d)$ is a subcomplex of the full de Rham complex,
\[\xymatrix{
0\ar[r] &K^1\ar[r]^d\ar[d] &K^2\ar[r]^d \ar[d]&\cdots \\
{\oplus_{i=0}^{2n}}\w^i \ar[r]^d& \oplus_{i=0}^{2n}\w^i \ar[r]^d&\oplus_{i=0}^{2n}\w^i
\ar[r]&\cdots
}
\]
We denote  by $H^i(K^*)$ the cohomology group of the complex $(K^*, d)$.
It follows from the Hodge decomposition and the Lefschetz decomposition that 
the map $p^i : H^i(K^*) \to \oplus_{j=0}^{2n}H^j_{dR}(X)$ is injective for $i=1,2$.
\\
Let $(de^{z(t)}\psi)_{[k]}$ denotes the term of $de^{z(t)}\psi$ of degree $k$ in $t$. The first term is given by 
$$
(d e^{z(t)}\psi )_{[1]} = da\psi + d b_1\psi,
$$
where $da\psi =(d( \b+\ol\b)\ome^2)\w\psi \in \(\w^{2,1}\oplus\w^{2,1}\)\w\psi$.
Thus $da\psi \in K^2$ is $d$-exact. Since the map $p^2$ is injective, 
the class $[da\psi]\in H^2(K^*)$ vanishes and we have 
a solution $b_1\in K^1$ of the first equation 
$da\psi + d b_1\psi=0$. 
\\
Next we consider an operator $e^{-z(t)} d e^{z(t)}$ acting on differential forms, where 
$z(t)=\log e^{at}e^{b(t)}$.
It follows that the operator $e^{-z(t)} d e^{z(t)}$ is a Clifford-Lie operator of order $3$ which is locally written in terms of the Clifford algebra valued Lie derivative, 
\bgn{equation}\label{eq: 2.6}
e^{-z(t)} d e^{z(t)}= \sum_i E_i\L_{v_i}+ N_i,
\end{equation}
where $\L_{v_i}$ denotes the Lie derivative by a vector $v_i$ and $E_i \in \CL^1, \, \,
N_i\in \CL^3$ ({\it cf }definition 2.2 in \cite{Go2}).
We find an open covering $\{ U_\a\}$ of $X$ with  a non-vanishing holomorphic $n$-form 
$\Ome_\a$ on each $U_\a$. 
We denote by  $\Phi_\a$ the pair $(\Ome_\a, \psi)$. 
Since the set of almost \K\"ahler structures is invariant under the action of diffeomorphisms,
the Lie derivative of $\Phi_a$ by a vector field $v$ is given by
$$
\L_{v}\Phi_\a = a_\a\cdot\Phi_\a =(a_\a\cdot\Ome_\a, a_\a\cdot\psi),
$$
for a section $a_\a\in \CL^2$ on $U_\a$.
It follows from (\ref{eq: 2.6}) that there is a $h_\a\in \CL^3$ such that 
$$
e^{-z(t)} d e^{z(t)}\cdot\Phi_\a=h_\a\cdot\Phi_\a=(h_\a\cdot\Ome_\a, h_\a\cdot\psi).
$$
Since $b(t)\in \wtil\ker^1$ and $\Ad_{e^{z(t)}}\J_J=\J_{\b t}$ is integrable, 
we have  
$$
de^{z(t)}\Ome_\a= E_\a\cdot e^{z(t)}\Ome_\a,
$$
for a $E_\a\in \CL^1$, which is the integrablity condition of $\J_{\b t}$ in terms of pure spinors.
Thus 
$$e^{-z(t)} d e^{z(t)}\cdot\Ome_\a= \wtil E_\a\cdot\Ome_a,$$ where 
$\wtil E_a=e^{-z(t)}E_a e^{z(t)}\in \CL^1$.
It follows from  $h_a\cdot\Ome_\a=\til E_\a\cdot\Ome_\a$ that 
$h_\a \in \wtil \ker^2$.
It implies that $h_\a\cdot\psi \in K^2$.

We shall find a solution $b(t)$ of the equation $(d e^{z(t)}\psi )=0$ by the induction on the degree $k$.
We assume that there exists a solution $b_j\in\ker^1$ for $0\leq j<k$ of the equation
$(d e^{z(t)}\psi )_{[i]}=0$, for all $0\leq i<k$.
Then we have 
\bgn{align}\(e^{-z(t)} d e^{z(t)}\psi\)_{[k]}=&\sum_{\stackrel {i+j=k}{\ss 0\leq i, j\leq k}} 
\(e^{-z(t)}\)_{[j]}\(d\, e^{z(t)}\psi\)_{[i]}\\
=&\(d\, e^{z(t)}\psi\)_{[k]}.
\end{align}
Since $\(e^{-z(t)} d e^{z(t)}\psi\)|_{U_\a}=\,h_\a\cdot\psi|_{U_\a}\in K^2$ for $h_\a\in \wtil\ker^2$ on each $U_\a$, it follows that 
$\(d\, e^{z(t)}\psi\)_{[k]}=\(h_\a\cdot\psi\)_{[k]} \in K^2$.
The $d$-exact form $\(d\, e^{z(t)}\psi\)_{[k]}$ is written as 
$$
\(d\, e^{z(t)}\psi\)_{[k]}=\frac1{k!}(db_k\cdot\psi)+\Ob_k,
$$
where $\Ob_k$ is also a $d$-exact form in $K^2$ which defined in terms of 
$a$ and $ b_j$ for $1\leq j <k$.
Since the map $p^2$ is injective, 
 it follows that the class $[\Ob_k]\in H^2(K^*)$ vanishes and we have a solution $b_k$ of the equation 
$\(d\, e^{z(t)}\psi\)_{[k]}=0$. 
By our assumption of the induction, we have a solution $b(t)$ in the form of formal power series, which can be shown to be a convergent series.\\
The cohomology group $H^1(K^*)$ is given by $H^{1,1}(X)$. Then
by applying theorem 3.2 in \cite{Go3}, we obtain a $2$-parameter family of deformations of \K\"ahler structures $(\J_{\b t}, \psi_{t,s})$,
\bgn{theorem}
There exists a family of solutions $b_s(t)$ parameterized by $s\in H^{1,1}(M)$, which gives rise to 
deformations of \K\"ahler structures $(\J_{\b t}, \psi_{t,s})$.
\end{theorem}
}

\large{
\section{Deformations of bi-Hermitian structures}
Let $(X, \ome)$ be an $n$ dimensional compact K\"ahler manifold with a Poisson structure $\b$ and a complex structure $J$.
Then we have deformations of \K\"ahler structures $\{\J_{\b t}, \psi_t\}$ as in section 2. According to theorem by Gualtieri, there is the one to one correspondence between \K\"ahler structures and bi-Hermitian structures with a torsion condition. 
A bi-Hermitian structure is a triple $( I^+, I^-, h)$ consisting of two complex structures $I^+$ and $I^-$ and a Hermitian structure $h$ with respect to both 
$I^+$ and $I^-$. Let $\ome^\pm$ be the Hermitian $2$-form and $\ol\pa^\pm$ the $\ol\pa$-operator with respect to 
$I^\pm$ respectively. Then the torsion condition is given by 
$$
d^+_c\ome^+=-d^-_c\ome^-=H,
$$
where $H$ is a $d$-exact $3$-form and $d_c^\pm =\sqrt{-1}(\pa -\ol\pa)$.
Let $z(t)$ be a solution of the equation $de^{z(t)}\psi =0$ as in section $2$, 
which gives rise to deformations of \K\"ahler structures $\{\J_{\b t}, \psi_t\}$,
where $e^{z(t)}=e^{at}\,e^{b(t)}$ and $a=\b+\ol \b$. 
Then we have the corresponding deformations of bi-Hermitian structures 
$\{( I^+(t), I^-(t), h_t)\}$, where 
$I^+(0)=I^-(0)=J$.
Let $b_1$ be the first term of power series $b(t)$ in $t$. 
On an open set $U$, we find a basis $\{ Z_i\}_{i=1}^n$ of vector fields of type $(1,0)$ with respect to the complex structure $J$. 
We denote by $\ol\theta^i$ the $1$-form of type $(0,1)$ defined by the interior product of 
$-\sqrt{-1}\ome$ by $Z_i$, 
$$
\ol\theta^i =-\sqrt{-1}\,\,i_{Z_i} \,\ome.
$$
We define $\ol E^\pm_i$ to be $Z_i\pm \ol\theta^i\in (\TT)\otimes\C$.
Then $b_1\in \CL^2$ acts on $\ol E^\pm_i$ by the adjoint action, 
$$
[b_1, \ol E^\pm_i]\in (\TT)\otimes\C.
$$
We denote by $\ol\b(\ol\theta^i)$ the vector filed given by the contraction of $2$-vector $\ol\b$ by $1$-from $\ol\theta^i$.
Then we have a deformed basis $\{Z_i^\pm(t)\}$ of vectors of type $(1,0)$ with respect to $I^\pm(t)$ on $U$ which is written by the followings up to degree $1$ in $t$,
\bgn{align}
&Z_i^+(t) \equiv Z_i +\(\ol\b(\ol\theta^i)+ \pi_{\ss TX}[b_1, \ol E^+_i]\) t,
\qquad (\text{\rm mod } t^2 )\\
&Z_i^-(t) \equiv Z_i +\(-\ol\b(\ol\theta^i)+ \pi_{\ss TX}[b_1, \ol E^-_i]\) t,
\qquad (\text{\rm mod } t^2 ),
\end{align}
where $\pi_{\ss TX} : \TT\to TX$ denotes the projection.
\bgn{lemma}For $a=\b+\ol\b$, there exists a solution $z(t)$ of the equation $de^{z(t)}\cdot\psi=0$ such that the first term
$b_1$ is a real $2$-form. 
\end{lemma}
\bgn{proof}
The first term of the equation $de^{z(t)}\psi=0$ is given by 
$$
da \cdot\psi + d b_1\cdot\psi =0.
$$
Then we have 
$$da\cdot\psi =d(\b +\ol \b)\cdot\psi =-\frac12(\b\cdot\ome^2 +\ol\b\cdot\ome^2)\psi,$$
where $\b\cdot\ome^2$ denotes the interior product of the $4$-form $\ome^2$ by the $2$-vector $\b$.
The $d$-exact form  $-\frac12d(\b\cdot\ome^2 +\ol\b\cdot\ome^2)$ is a real form of 
type $\w^{2,1}\oplus \w^{1,2}$. 
As in proposition 4.7 in \cite{Go1}, we have a real elliptic complex, 
\[
\xymatrix{
0\ar[r] & P^{1,1}_{\R} \ar[r]^d & \(\w^{2,1}\oplus\w^{1,2}\)_\R \ar[r] ^d
& \cdots,}
\]
whose cohomology groups are respectively given by 
the harmonic real primitive form $ \Bbb P^{1,1}_\R$ of type $(1,1)$ and 
the real part of the Dolbeault cohomology 
$\(H^{2,1}(X)\oplus H^{1,2}(X)\)_\R$. Hence we obtain a real $b_1\in  
P^{1,1}_\R$ with $da\cdot\psi +db_1\cdot\psi =0$.
Hence the result follows.
\end{proof}
Hence for $b_1\in \w^2T^*X$,
it follows from $\pi_{\ss TX}[b_1, \ol E^+_i]=0$ that 
the $Z^\pm_i(t)$ is given by 
\bgn{equation}\label{eq: KS-class}
Z^\pm_i(t)=Z_i\pm \ol\b\cdot\ol\theta^i t,\qquad (\text{mod } t^2).
\end{equation}
The contraction between $\b$ and $\sqrt{-1}\ome$ is written as 
$$
\sqrt{-1}\beta\cdot\ome=\sum_i (\b\cdot\theta)\ol\theta^i\in T^{1,0}\otimes\w^{0,1}.
$$
Since $\sqrt{-1}\b\cdot\ome$ is $\ol\pa$-closed, we have a class
 $\sqrt{-1}[\b\cdot\ome]\in H^{0,1}(X,T^{1,0})$.
Then it follows from (\ref{eq: KS-class}) that the infinitesimal tangent of deformations $\{I^+(t)\}$ and $\{I^-(t)\}$ are respectively  given by the classes $\sqrt{-1}[\b\cdot\ome]$ and $-\sqrt{-1}[\b\cdot\ome]\in H^{0,1}(X, T^{1,0})$. 
Hence we have 
\bgn{theorem}
Let $X$ be a compact K\"alher manifold with a Poisson structure $\b$.
The class $[\b\cdot\ome]\in H^{0,1}(X,T^{1,0})$ defined by
the contraction of $\b$ by a K\"ahler form $\ome$ gives rise to 
unobstructed deformations. In other words, we have a vanishing of the obstruction class,
$[\b\cdot\ome, \b\cdot\ome]=0\in H^{0,2}(X,T^{1,0})$.
\end{theorem}
\bgn{proof}
Let $\{\J_{\b t}, \psi_t\}$ be deformations of \K\"ahler structures as in section $2$ with the corresponding deformations of 
bi-Hermitian structures $\{I^+(t), I^-(t), h_t\}$. 
Then the class $\sqrt{-1}[\b \cdot\ome]\in H^{0,1}(X,T^{1,0})$ is 
the infinitesimal tangent of the deformations of $I^+(t)$. 
Hence we obtain a vanishing of the obstruction class,
$[\b\cdot\ome, \b\cdot\ome]=0\in H^{0,2}(X,T^{1,0})$.
\end{proof}
\bgn{example}{\rm
Let $M$ be a complex torus of $n$ dimension and $X$ the product of $M$ and 
the projective space $\C P^1$. Deformations of $X$ were explicitly studied in 
\cite{KS}. 
A holomorphic vector field on $\C P^1$ is written as a linear combination, 
$$ a\frac{\pa}{\pa \zeta}+b\zeta\frac{\pa}{\pa \zeta}
+c\zeta^2\frac{\pa}{\pa \zeta},
$$
where $\zeta$ is the affine coordinates of $\C P^1$ and $a,b,c$ are constants. 
Let $\{z_1, \cdots, z_n\}$ be the coordinates of complex torus $M$. 
Then every representative $p$ of $H^{0,1}(X,T^{1,0})$ is given in the from,
$$
p=\sum_i\(a_i\frac{\pa}{\pa \zeta}+b_i\zeta\frac{\pa}{\pa \zeta}
+c_i\zeta^2\frac{\pa}{\pa \zeta}\)d\ol z_i+
\sum_{j,k}\lam_{j k}\frac{\pa}{\pa z_j}d \ol z_k,
$$
where $\lam_{jk}$ are constants.
We define an $n\times 3$ matrix $P$ by
$$ P=
\bgn{pmatrix}
a_1&b_1&c_1\\
a_2&b_2&c_2\\
\vdots&\vdots&\vdots\\
a_n&b_n&c_n
\end{pmatrix}
$$
Then we see that the class of obstruction $[p, p]$ vanishes if and only if the rank of the matrix $P$ is less than or equal to $1$.
On the other hand, every holomorphic $2$-vector on $X$ is given by 
$$
\b=\sum_i\(a_i\frac{\pa}{\pa \zeta}+b_i\zeta\frac{\pa}{\pa \zeta}
+c_i\zeta^2\frac{\pa}{\pa \zeta}\)\w\frac{\pa}{\pa z_i}+
\sum_{j,k}\lam_{j k}\frac{\pa}{\pa z_j}\w\frac{\pa}{\pa z_k}.
$$
Then the Schouten bracket $[\b, \b]$ also vanishes if and only if the rank of $P$
is less than or equal to $1$. Let $\ome$ be the K\"ahler form $\ome_{FS}+
\ome_M$, where $\ome_{FS}$ denotes the Fubini-Study form of $\C P^1$ and $\ome_M$ is the standard K\"ahler form of $M$.
Then the contraction $\b\cdot \ome $ is the representative $p$ and we have a surjective map 
$$
H^0(X, \w^2T^{1,0})\to H^{0,1}(X, T^{1,0}).
$$
Let $\Lam=\{\ome_\a\}_{\a=1}^{2n}$ be the discrete lattice of maximal rank $2n$ in $\C^{n}$ with $M=\C^{n}/\Lam$, where 
$\ome_\a=(\ome_{\a1},\cdots \ome_{\a n})$.
We denote by $V_i$ the holomorphic vector field
$a_i\frac{\pa}{\pa \zeta}+b_i\zeta\frac{\pa}{\pa \zeta}+c_i\zeta^2\frac{\pa}{\pa \zeta}$.
Then $V_i$ generates the automorphism
$\exp{(V_i)}$ of $\C P^1$.
For each $\a$, we define an automorphism  $\rho_t(\ome_\a)$ by 
$$
\rho_t(\ome_\a)=\sum_i \exp{(\ome_{\a i}V_it)}
$$
In the case of  the rank of $P = 1$, $V_1,\cdots , V_n$ are commuting vector fields. Hence $\rho_t$ gives rise to a representation of 
$\Lam =\pi_1(M)$ on automorphisms of $\C P^1$.
By the family of representations $\{\rho_t\}$, we obtain deformations of complex fibre bundles $\{ X_t\}$ over the torus $M$ starting from the trivial bundle $X_0= M\times \C P^1$, 
$$
X_t = M\times_{\rho_t} \C P^1\to M.
$$
In \cite{KS}, It is shown that the infinitesimal deformations of $\{X_t\}$ is the class 
$[p]\in H^{0,1}(X,T^{1,0})$, where $\lam_{j k}=0$.
}
\end{example}
}

\bgn{thebibliography}{99}
\bibitem{A.G.G}
V.~Apostolov, P.~Gauduchon, G.~Grantcharov,
{\it Bihermitian structures on complex surfaces},
Pro. London Math. Soc. {\bf 79}(1999), 414-428,
Corrigendum: {\bf 92}(2006), 200-202
\bibitem{B.K}
S.~Barannikov and M.~Kontsevich,
{\it Frobenius manifolds and formality of Lie algebras of
polyvector fields}, International Math. Res. Notices, (4)(1998), 201-215
\bibitem{B.G.C}
H.~Bursztyn, M~Gualtieri and G.~Cavalcanti,
{\it Reduction of Courant algebroids and generalized complex }structures,
Adv. Math., 211 (2), (2007), 726--765,
Math.DG/0509640
\bibitem{BS}
  J. Barton and M.~Sti{\'e}non
{\it Generalized complex submanifols}
Math.DG/0603480
\bibitem{Ca}
G.~Cavalcanti,
{\it New aspects of $dd^c$-lemma},
Math.DG/0501406
\bibitem{Ch}
C.C.~Chevalley,
{\it The algebraic theory of Spinors}
Columbia University Press, 1954
\bibitem{Co}
T.~J.~Courant
{\it Dirac manifolds}
Trans.Amer.Math.Soc, {\bf 319}(2), (1990), 631-661
\bibitem{DZ}
J.~P.~Dufour and N.~T.~Zung 
{\it Poisson structures and their normal forms}
Progress in Math. Vol. 242, Birkh\"auser, 2000
\bibitem{F}
A.~Fujiki
,Bihermitain anti-self-dual structures on Inoue surfaces, 
preprint, (2007)
\bibitem{Go1}
R.~Goto,
{\it Moduli spaces of topological calibrations,
Calabi-Yau, hyperK\"ahler, G$_2$ and Spin$(7)$ structures},
International Journal of Math. {\bf 115}, No. 3(2004), 211-257
\bibitem{Go2}
R.~Goto,
{\it On deformations of generalized Calabi-Yau, hyperK\"ahler, G$_2$ and Spin$(7)$ structures},
Math.DG/0512211
\bibitem{Go3}
R.~Goto
{\it Deformations of \complex  and \K\"ahler structures}
Math. DG/0705.2495
\bibitem{Gu1}
M.~Gualtieri,
{\it Generalized complex geometry}
Math.DG/0703298
\bibitem{Gu2}
M.~Gualtieri,
{\it Hodge decomposition for generalized K\"ahler manifolds}, Math. DG/0409093
\bibitem{Gu3}
M.~Gualtieri,
{\it Branes and Poisson varieties}
Math.DG/0710.2719
\bibitem{Hi1}
N.~Hitchin,
{\it Generalized Calabi-Yau manifolds},
Q. J. Math.,{\bf  54}, 281-308(2003),
Math. DG/0401221
\bibitem{Hi2}
N.~Hitchin,
{\it  Instantons, Poisson structures and generalized K\"ahler geometry},
Commun. Math. Phys. 265(2006), 131-164
\bibitem{Hi3}
N.Hitchin,
{\it Bihermitian metrics on Del Pezzo surfaces},
Math.DG/060821
\bibitem{Hu}
D.~Huybrechts,
{\it Generalized Calabi-Yau structures, K3 surfaces and B-fields}, math.AG/0306132,
International Journal of Math.16(2005) 13
\bibitem{Ko}
K.~Kodaira
{\it Complex manifolds and deformations of complex structures}
Grundlehren der Mathematischen Wissenschaften, {\bf 283}, springer-Verlag,
(1986)
\bibitem{KS}
K.~Kodaira and D.C.~Spencer,
{\it  On deformations of complex, analytic structures I,II}
Ann. of Math., 67(1958), 328-466
\bibitem{KSIII}
K.~Kodaira and D.C.~Spencer,
{\it On deformations of complex analytic structure, III.
stability theorems for complex structures}
Ann. of Math., 71(1960),43-76
\bibitem{L.T}
Y.~Lin and S.~Tolman,
{\it Symmetries in generalized K\"ahler geometry},
Commun. Math. Phys (to appear), Math. DG/0509069
\bibitem{L.P}
Z.~J.~Liu and Ping.~Xu
{\it On quadratic Poisson structures}
Lett. Math. Phy.{\bf 26}, (1992), No. 1, 33-42
\bibitem{L.W.P}
Z.-J.~Liu, A.~Weinstein and Ping. ~Xu, 
{\it Manin triples for Lie bialgebroids},
J.Diff. Geom, {\bf 45},(1997) 547-574
\bibitem{Na}
Y.Namikawa
{\it Poisson deformations of affine symplectic varieties}
Math.AG/0609741
\bibitem{OM}
Oren Ben-Bassart and Mitya Boyarchenko
{\it Submanifolds of generalized complex manifolds},
 J. of Symplectic Geom., 2 (3) (2004), 309-355
\bibitem{Pol}
A.~Polishchuk
{Algebraic geometry of Poisson varieties}
J. Math.Sci, Vol.84, No. 5, 1997, 1413-1444
\bibitem{Sa}
F.~Sakai,
{\it Anti-Kodaira dimension of ruled surfaces},
Sci. Rep. Saitama Univ. {\bf 2}(1982) 1-7
\bibitem{Se}
J. P. Serre,
{\it Lie Algebra and Lie groups}
Lecture Notes in Mathematics 1500,
Springer-Verlag
\bibitem{Ti}
G.~Tian,
{\it Smoothness of the universal deformations spaces of compact Calabi-Yau manifolds and its Peterson-Weil metric},
Mathematical Aspect of string theory(ed. S.T.~Yau),(1987),
World Scientific Publishing co., Singapore, 629-646
\bibitem{V}
Izu.~Vaisman
{\it Reduction and submanifolds of generalized complex manifolds}
Diff. Geom. Appl., 25 (2007), 147-166.
\end{thebibliography}
\end{document}